\documentclass[11pt, leqno]{amsart}
\usepackage{amsthm, amsmath}
\usepackage{amssymb} 
\usepackage{amscd}
\usepackage[curve,matrix,arrow,frame,tips]{xy}
\usepackage{verbatim}

\newtheorem{thm}{Theorem}[section]
\newtheorem{prop}[thm]{Proposition}
\newtheorem{lemma}[thm]{Lemma}

\newtheorem{Remark}[thm]{Remark}

\newcounter{ex}[section]

\newcommand{\C}{{\bf C}}

\newcommand{\Q}{{\bf Q}}

\newcommand{\Z}{{\bf Z}}

\newcommand{\ti}{\tilde}

 \renewcommand{\O}{{\mathcal O}}

\newcommand{\Pic}{{\rm Pic}}

\def\thfill{\null\nobreak\hfill}

\def\endproof{\thfill\vbox{\hrule
  \hbox{\vrule\hbox to 5pt{\vbox to 5pt{\vfil}\hfil}\vrule}\hrule}}

\begin{document}

\title[]{Grothendieck-Riemann-Roch and the moduli of Enriques surfaces}
\author[G. Pappas]{G. Pappas}
\thanks{*Partially supported by  NSF Grant DMS05-01409.}
\address{Dept. of
Mathematics\\
Michigan State
University\\
E. Lansing\\
MI 48824-1027\\
USA}
\email{pappas@math.msu.edu}
 
\date{\today}
 
\maketitle 

\section{} A  (complex)  Enriques surface is
a projective smooth connected algebraic surface $Y$ over $\C$ with
${\rm H}^1(Y,\O_Y)={\rm H}^2(Y, \O_Y)=(0)$,  
$(\Omega^2_{Y})^{\otimes\, 2} \simeq \O_Y$  
but $\Omega^2_{Y} \not \simeq \O_Y$ ([CF]). 
In this   note we give  a short proof of the fact that the coarse moduli space of complex Enriques surfaces is quasi-affine. This was first shown by Borcherds [B] using the denominator function of a generalized Kac-Moody superalgebra (the  ``$\Phi$-function"). From this, it follows (see [BKPS])
that any complete family of complex Enriques surfaces is isotrivial.
Our observation here is that the  ingredient of Borcherds' theory 
of infinite products can  be completely removed from the proof and can be  replaced by a simple use of the Grothendieck-Riemann-Roch theorem. 

\medskip

\noindent {\bf Acknowledgment:} The author is most grateful to I. Dolgachev for 
his comments.

\section{}  We refer the reader to [CS] and [BPV] for background on Enriques and $K3$ surfaces.
Let $L=E_8(-1)\oplus E_8(-1) \oplus H\oplus H\oplus H$ be the $K3$ lattice with involution $\rho: L\to L$ given by 
$$
\rho((x, y, z_1, z_2, z_3))=(y, x, -z_1, z_3, z_2)
$$
Here $H=\Z\cdot e\oplus \Z\cdot f$ is the hyperbolic plane
lattice with $\langle e, f\rangle=1$, and $E_8$ denotes the positive definite 
orthogonal lattice corresponding to the Dynkin diagram of type $E_8$. 
(The parentheses $(m)$ mean that we modify the pairing by multiplying it
by the integer $m$.)
Let $L^+=\{l\in L\ |\ \rho(l)=l\}$, $L^{-}=\{l\in L\ |\ \rho(l)=-l\}$.
Then $L^+$ is isometric to $E_8(-2)\oplus H(2)$; the unimodular lattice 
$L^+(1/2)$ can be identified with the Enriques lattice and is isometric to
$E_8(-1)\oplus H$. 
The lattice $L^-$ has signature $(2, 10)$.
We set
$$
\Omega=\{\omega\in {\bf P}(L\otimes\C)\ |\ \langle \omega, \omega\rangle=0,\ \langle \omega, \bar\omega\rangle>0\},
$$
$$
\Omega^-=\{\omega\in {\bf P}(L^-\otimes\C)\ |\ \langle \omega, \omega\rangle=0,\ \langle \omega, \bar\omega\rangle>0\}\ .
$$
for the $20$-dimensional (resp. $10$-dimensional) hermitian symmetric space 
for $O(3,19)$ (resp. for $O(2,10)$).
($\Omega^-$ has two $2$ connected components.) 
Also let $O(L)^\rho=\{g\in O(L)\ |\ g\cdot \rho=\rho\cdot g\}$, 
$$
O(L)^-={\rm restr}_{L^-}\left(O(L)^\rho\right) \ 
$$
where by $O$ we denote the orthogonal group of 
automorphisms of the lattice that preserve the form.
By work of Horikawa, $O(L)^-=O(L^-)$.
For 
 a vector $d$ of norm $\langle d, d\rangle=-2$ in $L^-$ let
$H_d$ be the divisor of points of $\Omega^-$ represented by 
vectors orthogonal to $d$ and set $\Omega^-_0=\Omega^--\left(\cup_d H_d\right) $.
The group $O(L)^-=O(L^-)$ acts on $\Omega^-_0$.

By the work of Namikawa and especially Horikawa ([H], [BPV])) on the Torelli theorem
for Enriques surfaces, the (coarse) moduli variety
$D^0$ of complex Enriques surfaces can be identified with the quotient
\begin{equation}\label{quo}
D^0=\Omega^-_0/O(L^-)\, .
\end{equation}
(This has the structure of a quasi-projective variety over $\C$.)
There is no ``universal'' family of Enriques surfaces over $D^0$
but this can be remedied by passing to a suitable finite cover:

\begin{prop}
There is a normal subgroup $\Gamma\subset O(L^-)$ of finite index 
such that:

(i) The group $\Gamma$ acts freely on $\Omega^-_0$.
The corresponding quotient 
$$
D^0_\Gamma=\Omega^-_0/\Gamma
$$ is a smooth variety  which supports a smooth family of 
Enriques surfaces $f:  Y\to D^0_\Gamma$.

(ii) There is a family of $K3$-surfaces $g: X\to D_\Gamma^0$
with a free involution $\iota:  X\to  X$ which respects $g$ and is such that 
$ Y= X/\langle \iota\rangle$. Denote by $p: X\to Y$ the covering morphism.

(iii) The line bundle of automorphic forms of weight $1$ on $\Omega^-$
(with fiber over the point $\omega \in \Omega^-\subset {\bf P}(L^-\otimes\C)$
given by the corresponding line in $L^-\otimes\C$) descends to the line bundle $\omega_\Gamma=g_*(\Omega^2_{ X/D^0_\Gamma})$ on $D^0_\Gamma$.
\end{prop}

\begin{Proof}
The existence of such a subgroup $\Gamma$ appears to be well-known to the experts;
for the sake of completeness we sketch  one possible construction here.
Let us consider the vector $v=(0,0,0, e+f,e+f)$
in $L^+\subset L$; we have $\langle v, v\rangle=4$. Let us consider the  groups
$\Gamma_v=\{g\in O(L)\ |\   g(v)=v\}<O(L)$, 
$$
\Gamma^\rho_v=\{g\in O(L)\ |\ \rho\cdot g=g\cdot \rho, g(v)=v\}<O(L), \quad \Gamma_v^-={\rm restr}_{L^-}(\Gamma^\rho_v)<O(L^-)\, .
$$
(The structure of $\Gamma_v^-$ is explained in [StI, Prop. 2.7].)
For odd $n\geq 3$, we   set $\Gamma_n=\{g\in  O(L) |\ g\equiv {\rm Id}\ {\rm mod}\ nL\}$
and set $\Gamma_{v,n}:=\Gamma_{v}\cap \Gamma_n$, $\Gamma^\rho_{v, n}:=\Gamma^\rho_{v}\cap \Gamma_n$. Similarly, denote by $\Gamma_{v, n}^-$ the restriction 
${\rm restr}_{L^-}(\Gamma^\rho_{v,n})$; using [StI, Prop. 2.7] and the fact that $n$ is odd we see that
$$
\Gamma_{v, n}^-=\Gamma^-_v\cap\{ h\in O(L^-)\ |\ h\equiv {\rm Id}\ {\rm mod}\ n\}\, .
$$
This is a normal subgroup   of finite index in $\Gamma^-_v$. By a well-known lemma
of Serre, $\Gamma_n$ and all our groups with a subscript $n$  have only trivial torsion.
The group $\Gamma_{v, n}$ acts on the homogeneous space 
$\Omega_v=\{\omega\in {\bf P}((\C\cdot v)^\perp)\ |\ \langle \omega, \omega\rangle=0,\ \langle \omega, \bar\omega\rangle>0\}$ with trivial stabilizers. Using the period morphism and the Torelli theorem (e.g [BR], [Fr]), we can identify the smooth quasi-projective quotient
$\Omega_v/\Gamma_{v, n}$ with the {\sl fine} moduli space of polarized $K3$ surfaces 
with a polarization of degree $4$ and $n$-level structure. Now notice  that the subgroup  $\Gamma^\rho_{v, n}$ preserves   $\Omega^-\subset \Omega_v$
and its action factors through $\Gamma^\rho_{v,n}\to \Gamma^-_{v,n}$. We have a morphism
$$
\Omega^-/\Gamma^-_{v, n}\rightarrow \Omega_v/\Gamma_{v, n}\ .
$$
The group $\Gamma^-_{v, n}$ acts freely on $\Omega^-$.
The quotient  $\Omega^-/\Gamma^-_{v, n}$ is also a smooth quasi-projective variety and 
supports a   family $\ti g: \ti X\to \Omega^-/\Gamma^-_{v, n}$ of $K3$-surfaces (obtained by pulling back the universal family over $\Omega_v/\Gamma_{v, n}$). In addition, we see (using  the  Torelli theorem again) that the family supports an involution $\ti\iota: \ti X\to \ti X$. Now restrict the family $\ti g$ over the complement of the images of the divisors $H_d$, i.e over
$$
D^0_{\Gamma^-_{v, n}}=\left(\Omega^--\left(\cup_d H_d\right)\right)/\Gamma^-_{v, n}.
$$
Using results of Nikulin  we see 
as in [Y, \S 1] that the involution $\ti \iota$ is fixed-point-free
over this complement. 
Denote by $g: X\to D^0_{\Gamma^-_{v, n}}$ the corresponding 
(restricted) family of $K3$'s with involution $\iota: X\to X$.  
The family $f: X/\langle \iota\rangle \to D^0_{\Gamma^-_{v, n}}$
is a family of Enriques surfaces with the desired properties.
For any subgroup $\Gamma\subset \Gamma^-_{v,n}$ of finite index the families $f$ and $g$ pull back to families over $D^0_\Gamma$ with the same properties; this allows
to make sure that $\Gamma$ can be taken to be normal in $O(L^-)$.\endproof
\end{Proof}
\medskip
 
For simplicity of notation we will set $S=D^0$, $T=D^0_\Gamma$
and $\omega_\Gamma=\omega$. 

The natural map $\pi: T\to S$ is a finite morphism of quasi-projective varieties
and $S$ is identified with the quotient $T/G$ where $G:=O(L^-)/\Gamma$. Now consider the line bundle $\omega ^{\otimes 2}$ over $T$. 

\begin{lemma} We have $f^*(\omega ^{\otimes 2})\simeq (\Omega^2_{ Y/T})^{\otimes 2}$ as line bundles on $Y$.
\end{lemma}

\begin{Proof}
Since $f: Y\to T$ is a family of Enriques surfaces, the line bundle $(\Omega^2_{ Y/T})^{\otimes 2}$ is trivial along the fibers of $f$.
Therefore, there is a line bundle $\delta$ on $T$ such that $(\Omega^2_{ Y/T})^{\otimes 2}\simeq f^*(\delta)$ and we have to show that $\delta\simeq \omega^{\otimes\, 2}$. Notice that 
$$
g^*\delta\simeq p^*f^*\delta\simeq p^*((\Omega^2_{ Y/T})^{\otimes 2})\simeq 
(\Omega^2_{ X/T})^{\otimes 2} \simeq  g^*\omega^{\otimes\, 2}\ .
$$
  The result now follows by applying ${\rm R}^0 g _*$ to both sides 
of $g^*\delta\simeq g^*\omega^{\otimes\, 2}$ after using the projection formula and ${\rm R}^0g_*(\O_X)\simeq\O_S$.\endproof  
\end{Proof}
\bigskip

Now let us apply the
Grothendieck-Riemann-Roch theorem ([Fu, Theorem 15.2]) to the proper morphism $f: Y\to T$
and to the trivial line bundle $\O_Y$  on $ Y$.
This gives
\begin{equation}\label{GRR}
{\rm ch}({\rm R} f_*\O_Y)=f_*({\rm ch}(\O_Y)\cdot {\rm Td}(T_f))=f_*({\rm Td}(T_f))
\end{equation}
in the Chow ring with rational coefficients ${\rm CH}^*(T)_\Q:={\rm CH}^*(T)\otimes_\Z\Q$.
In this formula, ${\rm ch}$ is the Chern character and ${\rm Td}$  the Todd class. Also ${\rm R} f_*\O_Y$ denotes the (virtual) total cohomology bundle,
and $T_f=(\Omega^1_{Y/T})^\vee$ is the relative tangent bundle of the smooth morphism $f$. On the right hand side of the formula, $f_*: {\rm CH}^{*}(Y)_\Q\to {\rm CH}^{*-2}(T)_\Q$ is the push forward homomorphism on Chow groups.
Set $v_1$ for the  class of $\omega $ in 
${\Pic}(T)$. By the above lemma, we have $f^*(v_1)=-2\cdot c_1(T_f)$.
For simplicity, let us denote by $c_i=c_i(T_{f})$ the Chern classes of 
the relative tangent bundle.
By reading off the degree $1$ component of the Grothendieck-Riemann-Roch identity 
(\ref{GRR}) we find
\begin{equation}
c_1({\rm R}  f_*(\O_Y))=\frac{1}{24}\cdot f_*\left(c_1\cdot c_2\right)
\end{equation} 
in ${\rm Pic}(T)\otimes\Q$.  
Since for the family of Enriques surfaces $f: Y\to T$, we have 
${\rm R}^0f_*(\O_Y)=\O_T$, ${\rm R}^if_*(\O_Y)=(0)$ if $i>0$,
we have ${\rm R}  f_*(\O_Y)=\O_T$. Hence, the above gives 
\begin{equation}\label{id1}
f_*(c_1\cdot c_2)=0
\end{equation}
in ${\rm Pic}(T)\otimes_\Z\Q$. Hence, we obtain 
$$
0=2\cdot f_*(c_1\cdot c_2)
=-f_*(f_*(v_1)\cdot c_2)=-v_1\cdot f_*(c_2)\, .
$$
By Noether's fomula $f_*(c_2)=12$ in the Chow group of codimension $0$ cycles 
${\rm CH}^0(T)$ and so
the class $v_1=c_1(\omega)$ is trivial in ${\rm Pic}(T)/{\rm torsion}\subset {\rm Pic}(T)\otimes_\Z\Q$. Therefore, the line bundle $\omega$ over $T$
is torsion. 

Now notice that automorphic forms of weight 
$n$ on $D^0_\Gamma$ are by construction sections of 
$$
(\pi_*\omega_{\Gamma}^{\otimes\, n})^G
$$
over $S=D^0$ (i.e  
$G=O(L^-)/\Gamma$-invariant sections of $\omega^{\otimes\, n}=\omega^{\otimes\, n}_{\Gamma}$
over $T=D^0_\Gamma$). Since the line bundle $\omega$ is torsion and the group $G$
is finite,  for $n$ sufficiently large and divisible, the line bundle  $\omega_{\Gamma}^{\otimes\, n}$ 
has a nowhere vanishing $G$-invariant global section. On the other hand, by Baily-Borel ([BB]), for all sufficiently large and divisible $n$ the line bundle of  
automorphic forms of weight 
$n$ on $D^0$ is very ample. These two facts together imply that
the moduli space $S=D^0$ is quasi-affine exactly as in [B].
As in loc. cit. it follows (using results of Sterk) that $D^0$ is 
in fact an affine variety with a copy of the affine line $\C$ removed.

\begin{Remark} 
{\rm
The above argument applies to any family of (classical) Enriques surfaces 
$f: Y\to T$, even if $T$ is a variety over a field of positive characteristic.
The conclusion is that  ${\rm R}^0f_*((\Omega^2_{Y/T})^{\otimes\, 2})$
is a torsion line bundle on $T$.   }
\end{Remark}

\end{document}